\documentclass[english]{article}
\usepackage[T1]{fontenc}
\usepackage[latin9]{inputenc}
\usepackage{amsmath}
\usepackage{amssymb}

\makeatletter
\newcommand{\lyxaddress}[1]{
\par {\raggedright #1
\vspace{1.4em}
\noindent\par}
}

\makeatother

\usepackage{babel}
\begin{document}

\title{A remark on the paper ``Redundancy of multiset topological spaces''}

\author{Moumita Chiney$^{1}$, S. K. Samanta$^{2}$}

\maketitle

\lyxaddress{Department of Mathematics, Visva-Bharati, Santiniketan-731235, West
Bengal, India.
e-mail: $^{1}$ moumi.chiney@gmail.com,$^{2}$ syamal\_123@yahoo.co.in}
\begin{abstract}

In this note the authors have raised the question regarding the validity of the main result in \cite{key-1} by setting an example.
\end{abstract}
In \cite{key-1} the author claimed that the multiset topology defined
on a set $X$ is equivalent to the general topology defined on the
product $X\times\mathbb{N}.$ But there are some doubts in this statement.
In fact, the main result in \cite{key-1} is the Theorem 3.4 and by
the following example it is seen that in proof of Theorem 3.4 the
results: \\
$(1)$ $\varphi(\sqcup\underset{i\in J}{M_{i}})=\underset{i\in J}{\cup}\varphi(M_{i}),$\\
$(2)$ $\varphi(\sqcup\underset{i\in J}{M_{i}})=\underset{i\in J}{\cup}\varphi(M_{i})$
and \\
$(3)$ $\varphi(M^{\Delta})=(\varphi(M))^{c}$ \\
do not hold. $\newline$ \\
$\mathbf{Example\:1.}$ Let $X=\{x,y,z\},\omega=4$ and $U,V\in[X]^{\omega}$
such that\\
$U=\{4/x,3/y,2/z\},$ $\: V=\{4/x,4/y,4/z\}.$ \\
Let $M_{1},M_{2}\in P^{*}(U)$ so that $M_{1}=\{4/x,3/y\}$ and $M_{2}=\{2/x,3/y\}.$
\\
Then $M_{1}\sqcup M_{2}=\{4/x,3/y\},\; M_{1}\sqcap M_{2}=\{2/x,3/y\}.$
The complement of $M_{2}$ in $[X]^{\omega}$ is $M_{2}^{\Delta}=\{2/x,1/y,4/z\}$
and the complement of $M_{2}$ with respect to $U$ is $\left(M_{2}^{\Delta}\right)_{U}=\{2/x,2/z\}.$\\
Therefore, $\varphi(M_{1})=\{(x,4),(y,3)\}$, $\varphi(M_{2})=\{(x,2),(y,3)\}$,
$\varphi(M_{1})\cup\varphi(M_{2})=\{(x,4),(x,2),(y,3)\}$, $\varphi(M_{1})\cap\varphi(M_{2})=\{(y,3)\},$
$\varphi(M_{1}\sqcup M_{2})=\{(x,4),(y,3)\}$ and $\varphi(M_{1}\sqcap M_{2})=\{(x,2),(y,3)\}.$
\\
Again, $\varphi(M_{2}^{\Delta})=\{(x,2),(y,1),(z,4)\}$ and $\varphi\left(\left(M_{2}^{\Delta}\right)_{U}\right)=\{(x,2),(z,2)\}.$\\
Now, $\varphi(U)=\{(x,4),(y,3),(z,2)\}$ and $\varphi(V)=\{(x,4),(y,4),(z,4)\}.$\\
Then $\varphi(U)\setminus\varphi(M_{2})=\{(x,4),(z,2)\},$ $\varphi(V)\setminus\varphi(M_{2})=\{(x,4),(y,4),(z,4)\},$\\
$X\times\mathbb{N}\setminus\varphi(M_{2})=X\times\mathbb{N}\setminus\{(x,2),(y,3)\}$
and\\
$\left(X\times\{1,2,3,4\}\right)\setminus\varphi(M_{2})=\{(x,1),(x,3),(x,4),(y,1),$
$(y,2),(y,4),(z,1),(z,2),$ \\
$(z,3),(z,4)\}.$$\newline$\\
From the above results it follows that none of $(1),(2)$ and $(3)$ in
Theorem 3.4 of \cite{key-1} holds.

\end{document}